\documentclass[final,onefignum,onetabnum]{siamart171218}

\usepackage{lipsum}
\usepackage{amsfonts}
\usepackage{graphicx}
\usepackage{epstopdf}
\usepackage{algorithmic}
\ifpdf
  \DeclareGraphicsExtensions{.eps,.pdf,.png,.jpg}
\else
  \DeclareGraphicsExtensions{.eps}
\fi

\usepackage{autonum}

\newcommand{\Euler}{\genfrac{\langle}{\rangle}{0pt}{}}
\newcommand{\E}{\mathrm{E}}

\newcommand{\Beta}{\mathrm{Beta}}
\renewcommand{\P}{\mathrm{P}}

\newcommand{\citep}{\cite}

\newsiamremark{remark}{Remark}
\newsiamremark{hypothesis}{Hypothesis}
\crefname{hypothesis}{Hypothesis}{Hypotheses}
\newsiamthm{claim}{Claim}

\headers{Column Sums for Binomial and Eulerian Numbers}{M. Li and R. Goldman}

\title{Limits of Sums for Binomial and Eulerian Numbers and their Associated Distributions}

\author{Meng Li\thanks{Department of Statistics, Rice University, Houston, TX
  (\email{meng@rice.edu}, \url{http://meng.rice.edu}).}
\and Ron Goldman\thanks{Department of Computer Science, Rice University, Houston, TX 
  (\email{rng@cs.rice.edu}, \url{https://www.cs.rice.edu/\~rng/}).}}

\usepackage{amsopn}

\usepackage{enumitem}
\makeatletter
\def\namedlabel#1#2{\begingroup
    #2%
    \def\@currentlabel{#2}%
    \phantomsection\label{#1}\endgroup
}
\makeatother

\ifpdf
\hypersetup{
  pdftitle={Column Sums for Binomial and Eulerian Numbers},
  pdfauthor={M. Li and R. Goldman}
}
\fi

\begin{document}

\maketitle

\begin{abstract}
We provide a unified, probabilistic approach using renewal theory to derive some novel limits of sums for the normalized binomial coefficients and for the normalized Eulerian numbers. We also investigate some corresponding results for their associated distributions -- the binomial distributions for the binomial coefficients and the Irwin-Hall distributions (uniform B-splines) for the Eulerian numbers.
\end{abstract}



\section{Introduction}\label{sec:introduction}
\begin{figure}[h!t!]
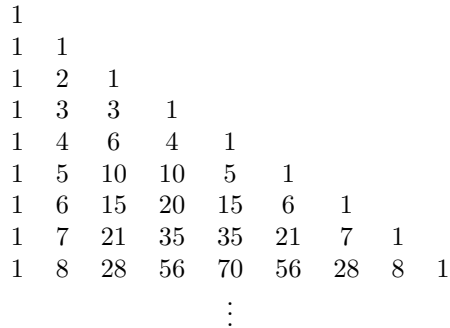

	\centering
	\begin{tabular}{ccccccccc}
		1&       &      &       &       &   &   &    &  \\
		1&      1&      &       &       &   &   &   &   \\
		1&      2&     1&       &       &   &   &   &   \\
		1&      3&     3&     1&        &   &   &   &   \\
		1&      4&     6&     4&     1  &   &   &   &   \\
		1&      5&    10&    10&    5   & 1 &   &   &   \\
		1&      6&    15&    20&   15   & 6 &   1&  &   \\
		1&      7&    21&    35&    35  & 21&  7 &  1&  \\
		1&      8&    28&    56&    70  & 56&  28&  8&  1  \\
		\multicolumn{9}{c}{$\vdots$}
	\end{tabular}
	\caption{Pascal's triangle -- levels 0 to 8 -- before normalization.}
	\label{fig:pascals triangle}
\end{figure}

Start with Pascal's triangle, the binomial coefficients ${n \choose k}$ arranged in a triangular array as in Figure~\ref{fig:pascals triangle}. Now normalize each row to sum to one: for each $n$ divide the $n$th row by $2^n$. After this normalization we can ask: what are the sums of the columns? The entries in the first column form the geometric progression $1, 1/2, 1/4, \ldots$ which sums to 2. In fact, it happens that when all the rows are normalized to sum to one, all the columns sum to two. But why 2? Equivalently we could ask: what are the sums of the long diagonals going from the upper left to the lower right? By symmetry the sums along these long diagonals are the same as the sums of the columns, so if we can sum the columns we can sum these long diagonals. But what about the short diagonals -- that is, what about the sums along the short diagonals that go from lower left to upper right? If we compute the sums along these diagonals, we see the Fibonacci series: $1, 1, 2, 3, 5, 8, 13, \ldots$. But these Fibonacci numbers emerge before we normalize the rows. We could also ask: what are these sums after we normalize the rows to sum to one? Then the series becomes: $1, 1/2, 3/4, 5/8, 11/16, 21/32,\ldots$. In fact, we shall show in Section~\ref{sec:binomial} that this series converges to 2/3. But why 2/3?

Consider next the Eulerian numbers $\Euler{n}{k}$ arranged again in a triangular array as in Figure~\ref{fig:eulers triangle}. For the Eulerian numbers we can also normalize each row to sum to one: for each $n$ divide the $n$th row by $n!$. Now once again we can ask: what are the sums of the columns? After normalization, the sum of the first column is the Taylor expansion of $e$, but the other columns certainly do not sum to $e$. What then can we say about the sums of these columns? By symmetry the sums along the long diagonals are the same as the sums of the columns, so if we can sum the columns we can sum these long diagonals. But again what about the short diagonals -- that is, what about the sums of the short diagonals that go from lower left to upper right after we normalize each row to sum to one? Are these sums in any way related to the corresponding sums in the normalized version of Pascal's triangle?

\begin{figure}[h!t!]
	\centering
	\begin{tabular}{ccccccccc}
		1 &&&&&&&&\\
		1   &      0&&&&&&&\\
		1   &      1   &      0&&&&&&\\
		1   &      4   &      1 &           0&&&&&\\
		1   &     11   &    11  &         1   &         0&&&&\\
		1   &     26   &    66  &        26  &         1  &         0&&&\\
		1   &     57   &   302  &      302  &       57 &         1 &      0&&\\
		1   &    120  &   1191  &   2416   &    1191 &     120 &   1  &    0&\\
		1   &    247  &   4293  &   15619  &  15619  &  4293 & 247 &  1  &   0  \\
		\multicolumn{9}{c}{$\vdots$}
	\end{tabular}
	\caption{The Eulerian numbers -- levels 0 to 8 -- before normalization.}
	\label{fig:eulers triangle}
\end{figure} 

The binomial coefficients ${n \choose k}$ count the number of subsets of $\{1,\ldots,n\}$ of exact order $k$, and satisfy the recurrence 
$$
{n \choose k} = {n - 1 \choose k} + {n - 1 \choose k - 1}.  
$$
In contrast, the Eulerian numbers $\Euler{n}{k}$ count the number of permutations of $\{1,\ldots,n\}$ with exactly $k$ ascents~\citep{grahamconcrete} and satisfy the recurrence
$$
\Euler{n}{k} = (k + 1) \Euler{n - 1}{k} + (n - k) \Euler{n - 1}{k - 1}. 
$$
(Here we follow the standard conventions that ${n \choose k}$ and $\Euler{n}{k}$ are both zero for $n < k$ and that ${0 \choose 0}$ and $\Euler{0}{0}$ are both equal to one.) At first glance then there is no reason to suspect any deep connection between these two triangular arrays of integers. But initial impressions can be deceiving; perhaps we can see some hidden connections if we look at some pictures.

\begin{figure}[h!t!]
	\centering
	\begin{tabular}{c @{\hskip 0.5in} c}
		\includegraphics[width = 0.3\linewidth]{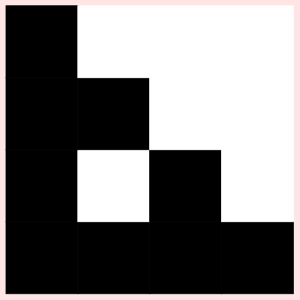}   & \includegraphics[width = 0.3\linewidth]{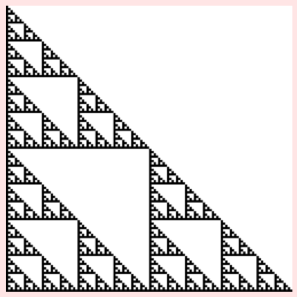} \\
		(a) Levels 0 to 3	& (b) Levels 0 to 127 \\
	\end{tabular}
	\caption{Pascal's triangle depicted by representing each odd number with a black square and each even number with a white square:  left -- levels 0 to 3, right -- levels 0-127.  As the number of levels increases, the Sierpinski triangle appears to emerge.}
	\label{fig:fractal.binomial}
	\vspace{0.2in}
	\begin{tabular}{c @{\hskip 0.5in} c}
		\includegraphics[width = 0.3\linewidth]{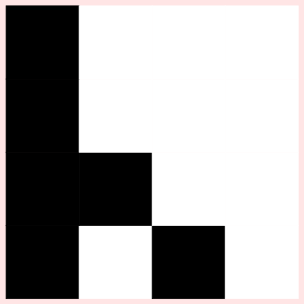}   & \includegraphics[width = 0.3\linewidth]{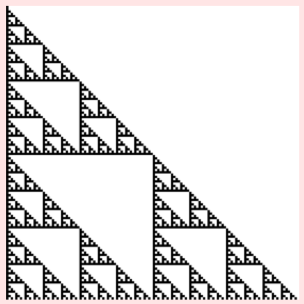} \\
		(a) Levels 0 to 3	& (b) Levels 0 to 127 \\
	\end{tabular}
	\caption{The Eulerian numbers depicted by representing each odd number with a black square and each even number with a white square: left -- levels 0 to 3, right -- levels 0-127.  Once again as the number of levels increases, the Sierpinski triangle appears to emerge.}
	\label{fig:fractal.eulerian}
\end{figure}

We can depict 2-dimensional arrays of integers in the following fashion: represent each odd integer by a black square and each even integer by a white square. Applying this approach to the binomial coefficients and to the Eulerian numbers generates Figures~\ref{fig:fractal.binomial} and~\ref{fig:fractal.eulerian}.

Recurrences often generate fractals. While level for level Figures~\ref{fig:fractal.binomial} and~\ref{fig:fractal.eulerian} are different -- compare, for example, levels 0 to 3 of the binomial coefficients with levels 0 to 3 of the Eulerian numbers -- in the large both arrays look very much the same: they both appear to incarnate the same fractal, the Sierpinski triangle. This long-term likeness suggests that in some limiting fashion these two arrays exhibit similar behaviors. \textit{The goal of this paper is to explore some limiting connections between sums of binomial coefficients and sums of Eulerian numbers along with corresponding results for their associated distributions: the binomial distributions for the binomial coefficients and the Irwin-Hall distributions (uniform B-splines) for the Eulerian numbers.} In particular, we are going to provide a unified, probabilistic approach using renewal theory to derive the following limiting identities: 

\begin{description}
	\item A. Binomial Coefficients
	\begin{description}
		\item[\namedlabel{eq:binomial.cs}{A1}] $\lim_{k \rightarrow \infty}  \sum_{n = 0}^{\infty }\frac{1}{2^n} {n \choose k} = 2$ \hfill (Columns)
		\item[\namedlabel{eq:binomial.sd}{A2}] $\lim_{n \rightarrow \infty} \sum_{k \geq 0} \frac{1}{2^{n - k}}{{n - k}\choose{k}} =  \frac{2}{3}$ \hfill (Short Diagonals)
		\item[\namedlabel{eq:binomial.as}{A3}] $\lim_{k \rightarrow \infty}  \sum_{n = 0}^{\infty} (-1)^n \frac{1}{2^n} {n \choose k} = 0$ \hfill (Alternating Sums)
	\end{description}
	\item B. Eulerian Numbers 
	\begin{description}
		\item[\namedlabel{eq:eulerian.cs}{B1}] $\lim_{k \rightarrow \infty} \sum_{n = 0}^{\infty} \frac{1}{n!} \Euler{n}{k} = 2$ \hfill (Columns)
		\item[\namedlabel{eq:eulerian.sd}{B2}] $\lim_{n \rightarrow \infty} \sum_{k \geq 0} \frac{1}{(n - k)!}\Euler{n - k}{k}  =  \frac{2}{3}$ \hfill (Short Diagonals)
		\item[\namedlabel{eq:eulerian.as}{B3}] $\lim_{k \rightarrow \infty}  \sum_{n = 0}^{\infty} (-1)^n \frac{1}{n!} \Euler{n}{k} = 0$ \hfill (Alternating Sums)
	\end{description}
	\item C. Bernstein Polynomials -- Binomial Distributions 
	\begin{description}
		\item[\namedlabel{eq:bernstein.cs}{C1}] $\lim_{k \rightarrow \infty} \sum_{n = 0}^{\infty} B^n_k(t) = \frac{1}{t}$ for $t \in (0, 1)$  \hfill (Columns)
		\item[\namedlabel{eq:bernstein.sd}{C2}] $\lim_{n \rightarrow \infty} \sum_{k \geq 0} B_k^{n - k}(t) = \frac{1}{1 + t}$ for $t \in (0, 1)$ \hfill (Short Diagonals) 
		\item[\namedlabel{eq:bernstein.as}{C3}] $\lim_{k \rightarrow \infty} \sum_{n = 0}^{\infty} (-1)^n B^n_k(t) = 0$ for $t \in (0, 1)$  \hfill (Alternating Sums)
	\end{description}
	\item D. $h$-Bernstein Polynomials -- P\'{o}lya-Eggenberger Distributions 
	\begin{description}
		\item[\namedlabel{eq:hbernstein.cs}{D1}] $\lim_{k \rightarrow \infty} \sum_{n = 0}^{\infty} B_k^n(t; h) =\frac{1 - h}{t - h}$ for $0 <  h < t < 1$ \hfill (Columns)
		\item[\namedlabel{eq:hbernstein.sd}{D2}] $\lim_{n \rightarrow \infty} \sum_{k \geq 0} B_k^{n - k}(t; h) = \int_0^1 \frac{x^{a - 1} (1 - x)^{b - 1}}{(1 + x) \mathrm{B}(a, b)} dx$ \hfill (Short Diagonals) \\
		for $t \in (0, 1)$ and $h > 0$, where $a = t/h, b = (1 - t)/h$
		\item D2a. $\lim_{n \rightarrow \infty} \sum_{k \geq 0} B_k^{n - k}(t; 1) = 2^{-t}$ for $t \in (0, 1)$ \hfill ($h$ = 1 in D2) 
		\item[\namedlabel{eq:hbernstein.as}{D3}] $\lim_{k \rightarrow \infty} \sum_{n = 0}^{\infty} (-1)^n B_k^n(t; h) = 0 $ for $0 <  h < t < 1$ \hfill (Alternating Sums)
	\end{description}
	\item E. Uniform B-splines -- Irwin-Hall Distributions 
	\begin{description}
		\item[\namedlabel{eq:bsplines.cs}{E1}] $\lim_{t \rightarrow \infty} \sum_{n =0}^{\infty}N_{0, n}(t) = 2$ \hfill (Columns)
		\item[\namedlabel{eq:bsplines.sd}{E2}] $\lim_{n \rightarrow \infty} \sum_{k \geq 0} N_{0, n - k}(k + t) = \frac{2}{3}$ for $ t > 0$ \hfill (Short Diagonals)
		\item[\namedlabel{eq:bsplines.as}{E3}]$\lim_{t \rightarrow \infty} \sum_{n =0}^{\infty} (-1)^n N_{0, n}(t) = 0$ \hfill (Alternating Sums)
	\end{description}
\end{description}

For the identities in B and E, we shall show that the convergence rate is polynomial with an arbitrarily large order, utilizing additional results from renewal theory. For the identities in A, C, and D, we shall show that analogous results hold for any fixed $k$ also invoking a probabilistic argument based on an infinite sequence of random variables, the same framework we use for the asymptotic identities. Below is a list of the non-asymptotic identities we will derive:  
\begin{description}
	\item A*. Binomial Coefficients 
	\begin{description}
		\item[\namedlabel{eq:binomial.cs.exact}{A1*}] $\sum_{n = 0}^{\infty }\frac{1}{2^n} {n \choose k} = 2$ \hfill (Columns)
		\item[\namedlabel{eq:binomial.sd.exact}{A2*}] $\sum_{k \geq 0} \frac{1}{2^{n - k}}{{n - k}\choose{k}} =  \frac{2}{3} + \frac{1}{3}\cdot \left(-\frac{1}{2}\right)^n$ \hfill (Short Diagonals)
		\item[\namedlabel{eq:binomial.as.exact}{A3*}] $\sum_{n = 0}^{\infty} (-1)^n \frac{1}{2^n} {n \choose k} = (-1)^k \frac{2}{3^{k + 1}}$ \hfill (Alternating Sums)
	\end{description}
	\item C*. Bernstein Polynomials -- Binomial Distributions 
	\begin{description}
		\item[\namedlabel{eq:bernstein.cs.exact}{C1*}] $\sum_{n = 0}^{\infty} B^n_k(t) = \frac{1}{t}$  \hfill (Columns)
		\item[\namedlabel{eq:bernstein.sd.exact}{C2*}] $\sum_{k \geq 0} B^{n-k}_k(t) = \frac{1 - (-t)^{n + 1}}{1 + t}$ \hfill (Short Diagonals) 
		\item[\namedlabel{eq:bernstein.as.exact}{C3*}]  $\sum_{n = 0}^{\infty} (-1)^n B^n_k(t) = \frac{(-t)^k}{(2 - t)^{k + 1}}$ \hfill (Alternating Sums)
	\end{description}
	\item D*. $h$-Bernstein Polynomials -- P\'{o}lya-Eggenberger Distributions 
	\begin{description}
		\item[\namedlabel{eq:hbernstein.cs.exact}{D1*}] $\sum_{n = 0}^{\infty} B_k^n(t; h) =\frac{1 - h}{t - h}$ for $0 < h < t < 1$ \hfill (Columns) 
		\item[\namedlabel{eq:hbernstein.sd.exact}{D2*}] $\sum_{k \geq 0} B_k^{n - k}(t; h) = \int_0^1 \frac{x^{a - 1} (1 - x)^{b - 1}}{(1 + x) \mathrm{B}(a, b)} dx +$ \hfill (Short Diagonals) \\
		\hspace*{0.3in}$(-1)^{n + 2} \int_0^1 \frac{x^{a + n} (1 - x)^{b - 1}}{(1 + x) \mathrm{B}(a, b)} dx$  \\
		for $t \in (0, 1)$ and $h >0,$ where $a = t/h, b = (1 - t)/h$
		\item[\namedlabel{eq:hbernstein.as.exact}{D3*}] $\sum_{n = 0}^{\infty} (-1)^n B_k^{n}(t; h) = (-1)^k \int_0^1 \frac{x^{a + k - 1} (1 - x)^{b - 1}}{(2 - x)^{k + 1} \mathrm{B}(a, b)} dx$ \hfill (Alternating Sums) \\
		for $0 < h < t < 1$, where $a = t/h, b = (1 - t)/h$
	\end{description}
\end{description}

Although it is reassuring that these identities involving binomial coefficients and Bernstein polynomials hold for any $k$ not only in the limit, the use of renewal theory provides an interesting link between binomial coefficients and Eulerian numbers and indicates the column sums in~\ref{eq:binomial.cs} and~\ref{eq:eulerian.cs} are expected to be 2, since 2 is the reciprocal mean of both the uniform distribution on $[0, 1]$ and the Bernoulli distribution with success probability 1/2. Likewise, the sums 2/3 of the short diagonals in~\ref{eq:binomial.sd} and~\ref{eq:eulerian.sd} are the reciprocal mean of a shifted version of these two distributions. The combinatorial interpretations of binomial coefficients and Eulerian numbers provide virtually no insight into revealing such a remarkable connection. Some identities such as~\ref{eq:eulerian.cs} seem straightforward using the general theory of renewal processes but are far from obvious or do not even appear promising otherwise (including the application of generating functions and enhanced approximation to the summand at each $n$).  

The rows of the sequences we shall study form distributions; the columns do not. While it may seem unnatural to study the columns rather than the rows, there has been some recent work to investigate similar column sums with good effect.~\cite{simsek2014generating} introduces generating functions for the columns of the binomial distribution and uses these generating functions to derive a variety of identities for the Bernstein polynomials, including formulas for sums, alternating sums, differentiation, degree elevation, and subdivision.~\cite{goldman2012generating} introduces generating functions for the columns of the uniform B-splines and then applies these generating functions to derive a collection of identities for uniform B-splines, including the Schoenberg identity, formulas for sums and alternating sums, for moments and reciprocal moments, for differentiation, for Laplace transforms, and for convolutions with monomials. Thus learning about the columns can also provide rich insights into the rows. It is partly in this spirit that we investigate the identities in this paper.

\section{Renewal theory}
Before we proceed with our proofs, we provide a brief review of renewal theory in stochastic processes~\citep[pp.\;358-373]{feller2008introduction}. A renewal process is a stochastic model for events that occur at random times. Let $X_1, X_2, \ldots$ be independent and identically distributed (IID) non-negative random variables following a distribution $F$ and let $S_n = \sum_{i = 1}^n X_i$ be their partial sum with $S_0 = 0$. We may interpret each $X_i$ as an \textit{interarrival time} and $S_n$ as the time of the $n$th arrival (or \textit{renewal}). 
The \textit{renewal measure} is defined by 
\begin{equation}
\label{eq:renewal.measure} 
U(A) = \sum_{n = 0}^{\infty} \P(S_n \in A), 
\end{equation}
for any $A$ that is a measurable subset on $[0, \infty)$. 
The renewal measure $U(A)$ is the expected number of arrivals in $A$ since 
\begin{align}
\label{eq:renewal.measure.expectation}
\sum_{n = 0}^{\infty} \P(S_n \in A)  & = \sum_{n = 0}^{\infty} \E\left(1\{S_n \in A\}\right) =  \E\left(\sum_{n = 0}^{\infty} 1\{S_n \in A\}\right) \\
& = \E\{ {\text{number of arrivals in }} A\},
\end{align}
where $1(\cdot)$ is the indicator function. 

Let $A = (x, x + \Delta]$ where $\Delta > 0$ is a fixed constant. The following renewal theorem, also known as Blackwell's theorem, states that the expected number of arrivals in $A$ is asymptotically proportional to the length of $A$ with proportionality constant $1/\mu$, where $\mu$ is the expectation of the random variables $X_i$. The precise statement of this theorem depends on whether or not the distribution $F$ is arithmetic: A distribution is called \textit{arithmetic} if it is supported on a set of the form $\{n\lambda: n \in \mathbb{N}\}$ for some $\lambda > 0$, and the largest such $\lambda$ is called the \textit{span} of the distribution.  
\begin{theorem}[Blackwell's renewal theorem]
	\label{thm:blackwell.non}
	If the distribution $F$ is not arithmetic, 
	$$\lim_{x \rightarrow \infty} U(x, x + \delta] = \delta/\mu$$ for any $\delta > 0$. 
	
	If the distribution $F$ is arithmetic with span $\lambda$, 
	$$\lim_{x \rightarrow \infty} U(x, x + \lambda] = \lambda/\mu.$$ 
\end{theorem}

As we shall see in the next section, a list of identities involving special numbers and special distributions can be derived using a unified probabilistic argument via renewal theory by specifying $A$ in Equation~\eqref{eq:renewal.measure} and the distribution $F$ of interarrival times. We consider mainly two classes of distributions: the uniform distribution supported on $[a, b]$ where $b > a$ and the Bernoulli distribution with success probability $p \in (0, 1)$. Normalized Eulerian numbers and uniform B-splines correspond to uniform distributions, while normalized Binomial coefficients, Bernstein polynomials, and $h$-Bernstein polynomials correspond to Bernoulli distributions. 

\newcommand{\Unif}{\mathrm{Uniform}}
\newcommand{\Bern}{\mathrm{Bernoulli}}
\newcommand{\Binom}{\mathrm{Binomial}}
\textbf{Notation.} We write $X \sim F$ if $X$ is a random variable following a distribution $F$. We use $\Unif(a, b)$ to denote the uniform distribution supported on $[a, b]$, whose probability density function is $f(x) = 1(x \in [a, b])$ and mean is $(a + b)/2$. We use $\Bern(p)$ to denote the Bernoulli distribution with success probability $p$, where the probability mass function is $\P(X = 1) = p$ (the trial succeeds) and $\P(X = 0) = 1 - p$ (the trial fails) and the mean is $p$. The sum of $n$ IID Bernoulli trials drawn from $\Bern(p)$ follows the binomial distribution $\Binom(n, p)$. 

We place a superscript on $X_i, S_n, U$ in a renewal process to emphasize their dependence on the interarrival distribution $F$. In particular, we use the superscript ``$[a, b]$" when $X_i \sim \Unif(a, b)$ and ``$(p)$" when $X_i \sim \Bern(p)$; we use  ``$(p) + 1$" when $X_i = X_i^* + 1$ where $X_i^* \sim \Bern(p)$, i.e., $X_i$ follows $\Bern(p)$ shifted by 1 satisfying $\P(X_i = 1) = 1 - p$ and $\P(X_i = 2) = p$. For example, the readers may see $X_i^{[0, 1]}$, $X_i^{[1, 2]}$, $X_i^{(1/2)}$ or $X_i^{(1/2) + 1}$, and similarly for $S_n$ and $U$. 

\section{Binomial coefficients and Eulerian numbers}
\subsection{Normalized binomial coefficients} 
\label{sec:binomial} 
The column sums of the normalized binomial coefficients are closely related to a renewal process due to the well known probabilistic interpretation of the normalized binomial coefficients using Bernoulli trials. Consider a renewal process in which the interarrival times $X_i \sim \Bern(1/2)$. Then 
\begin{equation}
\frac{1}{2^n} {n \choose k} = \P(S^{(1/2)}_n = k) = \P(S^{(1/2)}_n \in (k - 1, k]), 
\end{equation}
where $S_n^{(1/2)}$ is the time of the $n$th arrival and follows $\Binom(n, 1/2)$. 

Since the Bernoulli distribution is arithmetic with mean $\mu = 1/2$ and span $\lambda = 1$, Theorem~\ref{thm:blackwell.non} gives 
\begin{equation}
\lim_{k \rightarrow \infty}  \sum_{n = 0}^{\infty }\frac{1}{2^n} {n \choose k} =\lim_{k \rightarrow \infty} U^{(1/2)}(k - 1, k] = \frac{1}{1/2} = 2. \tag{A1}
\end{equation}

For the sums of the short diagonals, we consider a new renewal process with shifted Bernoulli interarrival times $X^{(1/2) + 1}_i$ that have the same distribution as $X^{(1/2)}_i + 1$, i.e., $X^{(1/2) + 1}_i = 2$ if the $i$th trial succeeds and $X^{(1/2) + 1}_i = 1$ if the $i$th trial fails. Short diagonals are closely related to this new renewal process because   
\begin{align}
\frac{1}{2^{n - k}} {{n - k}\choose{k}} & =\P(X^{(1/2)}_1 + \ldots + X^{(1/2)}_{n - k} \in (k - 1, k]) \\
& = \P((X^{(1/2)}_1 + 1) + \cdots + (X^{(1/2)}_{n - k} + 1) \in (n - 1, n]) \\
& = \P(X^{(1/2)+1}_1 + \cdots + X^{(1/2)+1}_{n - k}  \in (n - 1, n])  \\
& = \P(S^{(1/2)+1}_{n - k} \in (n - 1, n]). 
\end{align}
Therefore by Theorem~\ref{thm:blackwell.non}
\begin{align}
& \quad \; \lim_{n \rightarrow \infty} \sum_{k \geq 0} \frac{1}{2^{n - k}}{{n - k}\choose{k}} \\ & = \lim_{n \rightarrow \infty} \sum_{k \geq 0} \P(S^{(1/2)+1}_{n - k} \in (n - 1, n]) = \lim_{n \rightarrow \infty} \sum_{k \geq 0} \P(S^{(1/2)+1}_k \in (n - 1, n]) \label{eq:binomial.change.index} \\ & = \lim_{n \rightarrow \infty} U^{(1/2) +1}(n - 1, n] = \frac{1}{1 + \frac{1}{2}} = \frac{2}{3}, \tag{A2}
\end{align}
where the change of index in Equation~\eqref{eq:binomial.change.index} is guaranteed by the observation that $\P(S^{(1/2)+1}_k \in (n - 1, n]) = 0$ for $k \geq n$.

Equation~\eqref{eq:binomial.cs}, as well as its variants in A, C, and D, actually hold for any $k$. We provide a probabilistic proof of their non-asymptotic counterparts in Section~\ref{sec:non.asymptotic}. The use of renewal theory unites Binomial coefficients and Eulerian numbers under the same framework with two interarrival time distributions, and leads to a range of identities involving special distributions. In the next section, we shall elaborate on such connections. 

\subsection{Normalized Eulerian numbers} 
\label{sec:Eulerian}
Consider a renewal process with random interarrival time $X^{[0, 1]}_i \sim \Unif(0, 1)$. \cite{Tanny1973} provides a probabilistic interpretation for the normalized Eulerian numbers: for each integer $k$, 
\begin{equation} \label{eq:eulerian2unifom}
\frac{1}{n!} \Euler{n}{k} = \P(S^{[0, 1]}_n \in (k - 1, k]).
\end{equation} 
It follows by Equation~\eqref{eq:renewal.measure} that 
\begin{equation} \label{eq:Eulerian.renewal}
\sum_{n = 0}^{\infty} \frac{1}{n!} \Euler{n}{k} = \sum_{n = 0}^{\infty} \P(S^{[0, 1]}_n \in (k - 1, k]) = U^{[0, 1]}(k - 1, k]. 
\end{equation}
Equation~\eqref{eq:Eulerian.renewal} bridges the quantity $\sum_{n = 0}^{\infty} \frac{1}{n!} \Euler{n}{k}$ originating from Eulerian numbers with a renewal process, revealing the probabilistic interpretation of this column sum as the expected number of arrivals in the interval $(k - 1, k]$ when the interarrival time is uniformly distributed on $[0, 1]$. 

The uniform distribution on $[0, 1]$ is continuous, thus non-arithmetic, and has mean $\mu = 1/2$. Substituting $x = k$ and $\delta = 1$ in Theorem~\ref{thm:blackwell.non}, we obtain   
\begin{equation}
\lim_{k \rightarrow \infty} \sum_{n = 0}^{\infty} \frac{1}{n!} \Euler{n}{k} =  \frac{1}{1/2} = 2. \tag{B1}
\end{equation}

\textbf{Rate of convergence.} We have evaluated $\sum_{n = 0}^{\infty} \frac{1}{n!} \Euler{n}{k}$ numerically and observed that this sum converges very rapidly to 2. If $X \sim \Unif(0, 1)$, then $\E(X^{\alpha + 1}) = \int_0^1 x^{\alpha + 1} d x =  \frac{1}{\alpha + 2} < \infty$ for any $\alpha \geq 0$. According to Corollary 5.2 in~\cite{Konstantopoulos1999}, 
\begin{equation}
\sum_{n = 0}^{\infty} \frac{1}{n!} \Euler{n}{k} - 2 = o(k^{-\alpha}). 
\end{equation}
Therefore, the convergence rate in Equation~\eqref{eq:eulerian.cs} is polynomial with an arbitrarily large order. Below we calculate the difference $\sum_{n = 0}^{\infty} \frac{1}{n!} \Euler{n}{k} - 2$ for the first several $k$'s using ${\it Mathematica}$:


\begin{tabular}{ccccc} \hline 
 $k$  &  0 &  1 & 2& 3 \\
 $\sum_{n = 0}^{\infty} \frac{1}{n!} \Euler{n}{k} - 2$   & 0.71828 & $-4.8 \times 10^{-2}$ & $-4.2 \times 10^{-3}$ & $3.9 \times 10^{-5}$  \\ \hline
  $k$  & 4 & 5 & $\cdots$ & 50\\
 $\sum_{n = 0}^{\infty} \frac{1}{n!} \Euler{n}{k} - 2$   & $5.7 \times 10^{-5}$ & $5.1 \times 10^{-6}$ & $\cdots$ & $< 1.0 \times 10^{-45}$ \\ \hline 
\end{tabular}

\textbf{Short diagonals.} Short diagonals turn out to be related to a renewal process with interarrival times $X^{[1,2]}_i \sim \Unif(1, 2)$, which is identically distributed as $X^{[0, 1]}_i + 1$. The mean of $X^{[1,2]}_i  \sim \Unif(1, 2)$ is $\E X^{[1,2]}_i = \E X^{[0, 1]}_i + 1 = 3/2.$

Using Equation~\eqref{eq:eulerian2unifom} we find that 
\begin{align}
\frac{1}{(n - k)!}\Euler{n - k}{k}& =\P(X^{[0, 1]}_1 + \ldots + X^{[0, 1]}_{n - k} \in (k - 1, k]) \\
& = \P((X^{[0, 1]}_1 + 1) + \cdots + (X^{[0, 1]}_{n - k} + 1) \in (n - 1, n]) \\
& = \P(X^{[1, 2]}_1 + \cdots + X^{[1, 2]}_{n - k}  \in (n - 1, n])  = \P(S^{[1, 2]}_{n - k} \in (n - 1, n]). 
\end{align}
Therefore, 
\begin{align}
\sum_{k \geq 0} \frac{1}{(n - k)!}\Euler{n - k}{k}  & =\sum_{k \geq 0} \P(S^{[1, 2]}_{n - k} \in (n - 1, n]) \\ & = 
\sum_{k \geq 0} \P(S^{[1, 2]}_{k} \in (n - 1, n]) = U^{[1, 2]}(n-1, n], 
\end{align}
so by Theorem~\ref{thm:blackwell.non}
\begin{equation}
\lim_{n \rightarrow \infty} \sum_{k \geq 0} \frac{1}{(n - k)!}\Euler{n - k}{k} = \lim_{n \rightarrow \infty} U^{[1, 2]}(n-1, n] = \frac{1}{3/2} = \frac{2}{3},  \tag{B2}
\end{equation}
which is the same value as the analogous result for the binomial coefficients. 

\section{Bernstein polynomials and \textit{h}-Bernstein polynomials}
\subsection{Bernstein polynomials} \label{sec:bernstein} 
We consider a renewal process with random interarrival time $X^{(t)}_i \sim \Bern(t)$ where the success probability $t \in (0, 1)$. In view of the probabilistic interpretation of the Bernstein polynomials $B^n_k(t)$
\begin{equation}\label{eq:bernstein.renewal} 
B^n_k(t) =  {n \choose k} t^k (1 - t)^{n - k}  = \P(S_n^{(t)} = (k - 1, k]), 
\end{equation}
extensions of Equations~\eqref{eq:binomial.cs} and~\eqref{eq:binomial.sd} to Bernstein polynomials are immediately available. Since $\Bern(t)$ is arithmetic with span $\lambda = 1$ and mean $\mu = t$, a direct application of Theorem~\ref{thm:blackwell.non} gives 
\begin{equation}
\lim_{k \rightarrow \infty} \sum_{n = 0}^{\infty} B^n_k(t) = \lim_{k \rightarrow \infty} U^{(t)}(k - 1, k] = \frac{1}{t}. \tag{C1} 
\end{equation}

For the sums of the short diagonals in~\ref{eq:bernstein.sd}, we consider a new renewal process with shifted Bernoulli interarrival times $X^{(t) + 1}_i$ that have the same distribution as $X^{(t)}_i + 1$, i.e., $X^{(t) + 1}_i = 2$ if the $i$th trial succeeds and $X^{(t) + 1}_i = 1$ if the $i$th trial fails and has mean $\E X^{(t) + 1}_i = \E X^{(t)}_i + 1 = t + 1.$ The sums of the short diagonals follow from the same argument as used in the proof of~\ref{eq:binomial.sd} and~\ref{eq:eulerian.sd}: 
\begin{align}
B_k^{n - k}(t) & =\P(X^{(t)}_1 + \ldots + X^{(t)}_{n - k} \in (k - 1, k]) \\
& = \P((X^{(t)}_1 + 1) + \cdots + (X^{(t)}_{n - k} + 1) \in (n - 1, n]) \\
& = \P(X^{(t)+1}_1 + \cdots + X^{(t)+1}_{n - k}  \in (n - 1, n])  = \P(S^{(t)+1}_{n - k} \in (n - 1, n]), 
\end{align}
which by Theorem~\ref{thm:blackwell.non} gives 
\begin{align}
\lim_{n \rightarrow \infty} \sum_{k \geq 0} B_k^{n - k}(t) & = \lim_{n \rightarrow \infty} \sum_{k \geq 0} \P(S^{(t)+1}_{n -k} \in (n - 1, n]) \\
& = \lim_{n \rightarrow \infty} \sum_{k \geq 0} \P(S^{(t)+1}_k \in (n - 1, n]) = \lim_{n \rightarrow \infty} U^{(t) +1}(n-1, n] = \frac{1}{t + 1}, \tag{C2} 
\end{align}
noting that $\P(S^{(t)+1}_k \in (n - 1, n]) = 0$ for $k \geq n$.

\subsection{\textit{h}-Bernstein polynomials}\label{sec:h.bernstein}
The $h$-Bernstein polynomials are defined by 
\begin{equation}
B^n_k(t; h) = {n \choose k}\frac{\prod_{i = 0}^{k - 1} (t + i h) \prod_{i = 0}^{n - k - 1}(1 - t + i h)}{\prod_{i = 0}^{n - 1} (1 + i h)}, \quad k = 0, 1, 2, \ldots, n, 
\end{equation}
where $h \geq 0$. 
This formula is actually the probability density function of the P\'{o}lya Eggenberger distribution~\citep{eggenberger1923statistik,polya1930quelques}. This distribution reduces to the ordinary binomial distribution when $h = 0$, which has been discussed in the preceding section. Now we focus on positive $h$, when the P\'{o}lya Eggenberger distribution is a beta-binomial distribution with parameters $a = t/h$ and $b = (1 - t)/h$~\citep[Ch 6]{johnson2005univariate}. A beta-binomial distribution with parameters $(a, b)$ is the marginal distribution of $X$ if $(X | p) \sim \mathrm{Binomial}(n, p)$ and $p \sim \Beta(a, b)$, where $\Beta(a, b)$ is the Beta distribution having the probability density function $f(x) = \frac{x^{a - 1} (1 - x)^{b - 1}}{\mathrm{B}(a, b)}$ for $x \in (0, 1)$, and $\mathrm{B}(a, b)$ is the Beta function evaluated at $(a, b)$.

This interpretation of $B_k^n(t; h)$ using a mixture of binomial distributions means that $h$-Bernstein polynomials correspond to a renewal process with interarrival times $X^{(p)}_i$, where the success probability $p$ is a random draw from $\Beta(a, b)$. Therefore,  
\begin{equation}\label{eq:h.bernstein.renewal}
B_k^n(t; h) = \E_{p \sim \Beta(a, b) } \P(S^{(p)}_n \in (k - 1, k]) =: \int_0^1 \P(S^{(p)}_n \in (k - 1, k]) \frac{p^{a - 1} (1 - p)^{b - 1}}{\mathrm{B}(a, b)} dp.
\end{equation} 

Using the same argument as in the derivation for Bernstein polynomials but conditioning on the random success probability $p$, it follows that  
\begin{align}
\label{eq:exchange.limit}
\lim_{k \rightarrow \infty} \sum_{n = 0}^{\infty} B_k^n(t; h) & = \E_{p \sim \Beta(a, b) } \lim_{k \rightarrow \infty} \sum_{n = 0}^{\infty} \P(S^{(p)}_n \in (k - 1, k]) \\ & = \E_{p \sim \Beta(a, b) } \lim_{k \rightarrow \infty} U^{(p)}(k - 1, k] = \E_{p \sim \Beta(a, b) } \left(\frac{1}{p}\right) \\ & = \frac{a + b - 1}{a - 1} = \frac{1/h - 1}{t/h - 1} = \frac{1 - h}{t - h},
\tag{D1}
\end{align}
where the interchange of expectation and limit in Equation~\eqref{eq:exchange.limit} is guaranteed by the dominated convergence theorem~\citep[p.\;111]{feller2008introduction}. Here we require $a = t/h > 1$ to ensure that the expectation of $1/p$ exists, which means Equation~\eqref{eq:hbernstein.cs} holds for $0 < h < t < 1$. 

Asymptotic formulas for the sums of the short diagonals also hold for the $h$-Bernstein polynomials: 
\begin{align}
\lim_{n \rightarrow \infty} \sum_{k \geq 0} B_k^{n - k}(t; h) & = \E_{p \sim \Beta(a, b) }  \left[\lim_{n \rightarrow \infty} \sum_{n = 0}^{\infty} \P(S^{(p)+1}_n \in (n - 1, n])\right] \\
& = \E_{p \sim \Beta(a, b) } \left(\frac{1}{p + 1} \right) = \int_0^1 \frac{x^{a - 1} (1 - x)^{b - 1}}{(1 + x) \mathrm{B}(a, b)} dx, 
\label{eq:D2} \tag{D2}
\end{align}
where $0 < t < 1$, $h > 0$, $a = t/h$, $b = (1 - t)/h$. 

By Euler's integral representation of $_2F_1$~\citep[(1.6.1)]{koekoek2010hypergeometric}, the integral on the right hand side of Equation~\eqref{eq:D2} reduces to a hypergeometric function, so 
\begin{equation}
\lim_{n \rightarrow \infty} \sum_{k \geq 0} B_k^{n - k}(t; h) =\; _2F_1(1, t/h; 1/h; -1). 
\end{equation}
When $h = 1$, 
\begin{equation}
_2F_1(1, t/h; 1/h ; -1) =\; _2F_1(1, t; 1 ; -1) =\; _1F_0(t; ; -1) = 2^{-t},  \tag{D2a}
\end{equation}
for $0 < t < 1$. 

\section{B-splines} 
The application of renewal theory to Eulerian numbers can be extended to uniform B-splines~\citep{wang+:2010eulerian,he:2012eulerian}. Let $N_{0, n}(t)$ denote the B-spline of degree $n$ with knots at the integers $0, 1, \ldots, n+1$ and support $[0, n + 1]$. Let $\chi_{[0, 1]}(t)$ be the characteristic function over $[0, 1]$, i.e., $\chi_{[0, 1]}(t) = 1\{t \in [0, 1]\}$. Then $N_{0, n}(t)$ is the convolution of $\chi_{[0, 1]}(t)$ with itself $n + 1$ times. Since $\chi_{[0, 1]}(t)$ is the probability density function of $\Unif(0, 1)$, the convolution $N_{0, n}(t)$ is actually the probability density function of $S^{[0,1]}_{n + 1}$, or in other words, the probability density function of the Irwin-Hall distribution. Furthermore, 
\begin{equation}\label{eq:bsplines.prob} 
N_{0,n}(t) = \int_{-\infty}^{\infty} N_{0, n - 1}(x) \chi_{[0, 1]}(t - x) dx = \int_{t - 1}^{t} N_{0, n - 1}(x) dx =  \P(S^{[0, 1]}_{n} \in (t - 1, t]). 
\end{equation}
It follows from Equations~\eqref{eq:eulerian2unifom} and~\eqref{eq:bsplines.prob} that the normalized Eulerian numbers are the uniform B-splines evaluated at the integers. Moreover, in view of Theorem~\ref{thm:blackwell.non} and Equation~\eqref{eq:bsplines.prob}, it follows that 
\begin{equation}
\lim_{t \rightarrow \infty} \sum_{n =0}^{\infty}N_{0, n}(t) = 2. \tag{E1}
\end{equation}

Sums of short diagonals also have a counterpart for B-splines. Noting that
\begin{align}
& \quad \; N_{0, n - k}(k + t) =\P(X^{[0, 1]}_1 + \ldots + X^{[0, 1]}_{n - k} \in (k + t - 1, k + t]) \\
& = \P((X^{[0, 1]}_1 + 1) + \cdots + (X^{[0, 1]}_{n - k} + 1) \in (n + t - 1, n + t]) \\
& = \P(X^{[1, 2]}_1 + \cdots + X^{[1, 2]}_{n - k}  \in (n + t - 1, n + t])  = \P(S^{[1, 2]}_{n - k} \in (n + t - 1, n + t]), 
\end{align}
we have 
\begin{equation}
\lim_{n \rightarrow \infty} \sum_{k \geq 0} N_{0, n - k}(k + t) = \lim_{n \rightarrow \infty}U^{[1, 2]}(n + t - 1, n + t]  = \frac{1}{3/2} = \frac{2}{3}, \tag{E2} 
\end{equation}
for all $t$. 

\section{Contrasts and alternating sums} 
Theorem~\ref{thm:blackwell.non} also holds for a \textit{delayed renewal process}, where $S_0$ is a random variable other than a constant zero. This insight leads to an extension of Equation~\eqref{eq:eulerian.cs} to alternating sums, and more generally, \textit{contrasts}. We call a length $m$ vector $(c_0, c_1, \ldots, c_{m - 1})$ a {\it contrast} if 
\[
\sum_{i = 0}^{m-1} c_i = 0. 
\]
We use the same notation $c_k$ to denote its periodic extension, which is a sequence $\{c_k: k = 0, 1, 2, \ldots\}$ such that 
\[
c_k = c_i, \quad\text{if } k \equiv i \;(\bmod\; m). 
\]

For each $j = 1, 2, \ldots, m$, consider a delayed renewal process with delay $S_0^{(j)} = S^{[0,1]}_j$ and interarrival time $X_i^{(j)} = \sum_{i' = 0}^{m - 1} X^{[0,1]}_{m(i - 1) + i'  +  j + 1} = S^{[0,1]}_{mi + j} - S^{[0,1]}_{m(i - 1) + j}$ for $i = 1, 2, \ldots$. The corresponding partial sum $S_n^{(j)}$ satisfies  
\[
S_n^{(j)} = \sum_{i = 1}^n X_i^{(j)} + S^{[0,1]}_j = \sum_{i = 1}^n (S^{[0,1]}_{mi + j} - S^{[0,1]}_{m(i - 1) + j}) + S^{[0,1]}_j  = S^{[0,1]}_{mn + j}, 
\] 
and the expectation of interarrival times is $\E X_i^{(j)} = \sum_{i' = 0}^{m - 1} \E X^{[0,1]}_{m(i - 1) + i'  +  j + 1} = m/2$. 
A direct application of Theorem~\ref{thm:blackwell.non} to this delayed renewal process gives 
\begin{equation}
\lim_{k \rightarrow \infty} \sum_{n = 0}^{\infty} \P(S^{[0,1]}_{mn + j} \in (k - 1, k]) = \lim_{k \rightarrow \infty} \sum_{n = 0}^{\infty} \P(S_n^{(j)} \in (k - 1, k]) = \frac{1}{m/2} = \frac{2}{m}.
\end{equation}
Now 
\begin{equation}\label{eq:dummy12.1} 
\sum_{n = 0}^{\infty} c_n \frac{1}{n!} \Euler{n}{k} =\sum_{n = 0}^{\infty} c_n \P(S^{[0,1]}_n \in (k - 1, k])=  \sum_{j = 0}^{m - 1} c_j \sum_{n = 0}^{\infty} \P(S^{[0,1]}_{mn + j} \in (k - 1, k]). 
\end{equation}
Taking the limit of both sides in Equation~\eqref{eq:dummy12.1}, we obtain 
\begin{align}
\lim_{k \rightarrow \infty}  \sum_{n = 0}^{\infty} c_n \frac{1}{n!} \Euler{n}{k} & = \sum_{j = 0}^{m - 1} c_j \lim_{k \rightarrow \infty}  \sum_{n = 0}^{\infty} \P(S^{[0,1]}_{mn + j} \in (k - 1, k]) \\ & = \sum_{j = 0}^{m - 1} c_j \cdot \frac{2}{m} = \frac{2}{m} \sum_{j = 0}^{m - 1} c_j  = 0. 
\end{align}

A special case is alternating sums of the normalized Eulerian numbers when $m =2$ and $(c_0, c_1) = (1, -1)$, i.e., 
\begin{equation}
\lim_{k \rightarrow \infty}  \sum_{n = 0}^{\infty} (-1)^n \frac{1}{n!} \Euler{n}{k} = 0. \tag{A3}
\end{equation} 

Analogous results hold for normalized binomial coefficients following the same argument as above but replacing $\Unif(0, 1)$ by $\Bern(1/2)$. Thus  
\begin{equation}
\lim_{k \rightarrow \infty}  \sum_{n = 0}^{\infty} c_n \frac{1}{2^n} {n \choose k} = 0, 
\end{equation}
and 
\begin{equation} 
\lim_{k \rightarrow \infty}  \sum_{n = 0}^{\infty} (-1)^n \frac{1}{2^n} {n \choose k} = 0. \tag{B3}
\end{equation} 

Proofs of~\ref{eq:bernstein.as},~\ref{eq:hbernstein.as}, and~\ref{eq:bsplines.as} are similar to the proofs provided above either by varying the distribution of interarrival times (use $\Bern(t)$ where $t \in (0, 1)$ for~\ref{eq:bernstein.as} and a mixture of $\Bern(p)$ where $p \sim \Beta(a = t/h, b = (1 - t)/h)$ for~\ref{eq:hbernstein.as}) or by applying Theorem~\ref{thm:blackwell.non} to $(t - 1, t]$ for real-valued $t$ (for~\ref{eq:bsplines.as}). Thus we omit these proofs here.

\section{Non-asymptotic identities for normalized binomial coefficients and ($h$-)Bernstein polynomials}
\label{sec:non.asymptotic}
Here we derive non-asymptotic identities for fixed $k$ as a counterpart to the limits in A, C, and D, also through probabilistic arguments. 

The random variable $\sum_{n = 0}^{\infty} 1\{S^{(t)}_n = k\}$ is equal to one plus the number of trials that fail after $S_n$ first reaches $k$. Because the trails after $S_n$ are independent of $S_n$, $\sum_{n = 0}^{\infty} 1\{S^{(t)}_n = k\}$ has the same distribution as the number of trials before one success in an IID Bernoulli sequence, which is known to be a geometric distribution with parameter $t$ that has mean $1/t$. Consequently by Equation~\eqref{eq:renewal.measure.expectation}, 
\begin{equation}
\sum_{n = 0}^{\infty } {n \choose k} t^k (1 - t)^{n - k}  = \frac{1}{t}, \tag{C1*}
\end{equation}
for any nonnegative integer $k$. 
Substituting $t = 1/2$ into Equation~\eqref{eq:bernstein.cs.exact} leads to the column sum of the normalized binomial coefficients
\begin{equation}
\sum_{n = 0}^{\infty }\frac{1}{2^n} {n \choose k} = \frac{1}{1/2} = 2, \tag{A1*}
\end{equation}
for any nonnegative integer $k$. Equation~\eqref{eq:hbernstein.cs.exact} holds by using the equivalence between $h$-Bernstein polynomials and Beta-Binomial distributions as in Section~\ref{sec:h.bernstein}. 

For the sum of short diagonals, 
$$\sum_{k \geq 0} B_k^{n - k}(t) = \sum_{k \geq 0}  \P(S^{(t)+1}_{k} = n) = \E \left[ \sum_{k \geq 0} 1(S^{(t)+1}_{k} = n)\right] =: \E Z_n.$$
Since there exists at most one value of $k$ such that $S^{(t)+1}_{k} = n$ as $S^{(t)+1}_{k}$ is strictly increasing in $k$, the random variable $Z_n$ follows a Bernoulli distribution. Let the success probability be $a_n = \P(Z_n = 1)$. In the infinite sequence of binary trials $\{S^{(t)+1}_k\}$, conditioning on whether the event right before $\{Z_n = 1\}$ occurs is a success or a failure, we conclude that $a_n = a_{n - 1} (1 - t) + a_{n - 2} t$ for $n \geq 2$, where we let $a_0 = 1$. Solving this recurrence for $a_n$ gives 
\begin{equation}
\sum_{k \geq 0} B_k^{n - k}(t) = \frac{1 - (-t)^{n + 1}}{1 + t}. \tag{C2*}
\end{equation} 
Substituting $t = 1/2$ yields
\begin{equation} 
\sum_{k \geq 0} \frac{1}{2^{n - k}}{{n - k}\choose{k}}   = \frac{2}{3} + \frac{1}{3}\cdot \left(-\frac{1}{2}\right)^n. \tag{A2*}
\end{equation}
In view of $B_k^{n - k}(t; h) = \E_{p \sim \Beta(a, b) } B_k^{n - k}(p)$, which follows from  Equations~\eqref{eq:bernstein.renewal} and~\eqref{eq:h.bernstein.renewal} with $a = t/h$ and $b = (1 - t)/h$, Equation~\eqref{eq:bernstein.sd.exact} leads to 
\begin{align}
\sum_{k \geq 0} B_k^{n - k}(t;h) & = \E_{p \sim \Beta(a, b) } \sum_{k \geq 0} B_k^{n - k}(p) = \E_{p \sim \Beta(a, b) } \left[\frac{1 - (-p)^{n + 1}}{1 + p} \right] \\ &= \int_0^1 \frac{x^{a - 1} (1 - x)^{b - 1}}{(1 + x) \mathrm{B}(a, b)} dx + (-1)^{n + 2} \int_0^1 \frac{x^{a + n} (1 - x)^{b - 1}}{(1 + x) \mathrm{B}(a, b)} dx. \tag{D2*}
\end{align}

For the alternating sums, we just need to show the Bernstein polynomials in~\ref{eq:bernstein.as.exact}, then~\ref{eq:binomial.as.exact} and~\ref{eq:hbernstein.as.exact} will follow in the same way as we derive~\ref{eq:binomial.sd.exact} and~\ref{eq:hbernstein.sd.exact}. For Equation~\eqref{eq:bernstein.as.exact}, we first rewrite the left-hand side into an expectation 
\[
\sum_{n = 0}^{\infty} (-1)^n B_k^n(t) = \sum_{n = 0}^{\infty} (-1)^n \E 1\{S_n^{(t)} = k\} = \E \left[\sum_{n = 0}^{\infty} (-1)^n 1\{S_n^{(t)} = k\} \right] = \E Z_k. 
\] 
Consider an IID Bernoulli sequence with success probability $t$ and a derived random sequence $\{U_j: j = 1, \ldots\}$ where $U_j$ is the number of trials needed from the $(j - 1)$th success to the $j$th success. Then the $U_j$'s are IID following a geometric distribution with parameter $t$. Let $V_j$ be the number of the trial at which the $j$th success first occurs, i.e., $V_j = \sum_{j' = 1}^j U_{j'}$. Then, 
\begin{align}
Z_k & = \sum_{n = V_k}^{V_{k + 1} - 1} (-1)^n 1\{S_n^{(t)} = k\} =  \sum_{n = V_k}^{V_{k + 1} - 1} (-1)^n = (-1)^{V_k} \cdot 1\{ V_{k + 1} - V_k \text{ is odd}\} \\
& = \prod_{j = 1}^k (-1)^{U_j} \cdot 1\{ U_{k + 1} \text{ is odd}\}, 
\end{align}
which combined with the fact that the $U_j$'s are IID leads to  
\begin{equation}
\E Z_k = \prod_{j = 1}^{k} \E[(-1)^{U_j}] \P\{ U_{k + 1} \text{ is odd}\} = \{\E[(-1)^{U_1}]\}^k \P\{ U_{1} \text{ is odd}\}. 
\end{equation}
Let $A_1 = \{ U_1 \text{ is even}\}$ and $B_1 = \{ U_1 \text{ is odd}\} = A_1^c$. Then $\E[(-1)^{U_1}] = \P(A_1) - \P(B_1)$. Moreover it follows easily from the definition of $U_1$ that $\P(A_1) = \sum_{i = 0}^{\infty} (1 - t)^{2i + 1} t = t(1 - t)/(1 - (1 - t)^2) = (1 - t)/(2 - t)$ and $\P(B_1) = 1 - \P(A_1) = 1/(2 - t)$. Consequently, 
\begin{equation}
\E Z_k = [\P(A_1) - \P(B_1)]^k \P(B_1) = \frac{(-t)^k}{(2 - t)^{k + 1}}. \tag{C3*}
\end{equation}

%
%
%
%
%

%


\section*{Acknowledgments}
We thank Professor Plamen Simeonov for pointing out the connection between~\ref{eq:hbernstein.sd} and hypergeometric functions.

\bibliographystyle{siamplain}
\bibliography{library}
\end{document}